\documentclass[leqno]{article}
\usepackage{verbatim, amsmath, amscd, amssymb}
\usepackage[all]{xy}
\newcommand{\bU}{{\bf U}}

\newcommand{\cB}{{\mathcal B}}
\newcommand{\cC}{{\mathcal C}}
\newcommand{\cD}{{\mathcal D}}

\newcommand{\cL}{{\mathcal L}}

\newcommand{\cO}{{\mathcal O}}

\newcommand{\frg}{\mathfrak g}

\newcommand{\rC}{\mathrm{C}}

\newcommand{\rd}{\mathrm{d}}

\newcommand{\rs}{\mathrm{s}}

\newcommand{\bbE}{\mathbb E}

\newcommand{\bbN}{\mathbb N}

\newcommand{\bbR}{\mathbb R}

\newcommand{\bbZ}{\mathbb Z}

\newcommand{\Dist}{\mathrm{Dist}}

\newcommand{\gr}{\mathrm{gr}}

\newcommand{\hd}{\mathrm{hd}}

\newcommand{\im}{\mathrm{im}}

\newcommand{\ind}{\mathrm{ind}}

\newcommand{\op}{\mathrm{op}}

\newcommand{\pr}{\mathrm{pr}}\newcommand{\barpr}{\overline{\mathrm{pr}}}

\newcommand{\rad}{\mathrm{rad}}

\newcommand{\soc}{\mathrm{soc}}

\newcommand{\Mod}{\mathbf{Mod}}

\newcommand{\lbr}{\begin{bmatrix}}
\newcommand{\rbr}{\end{bmatrix}}
\newcommand{\for}{\bigcirc\kern-2.6ex \because}
\newcommand{\forb}{\bigcirc\kern-2.8ex \because}
\newcommand{\forbb}{\bigcirc\kern-3.0ex \because}
\newcommand{\forbbb}{\bigcirc\kern-3.1ex \because}

\newcommand\pf{\noindent {\bf Proof:  }}
\newtheorem{thm}{Theorem:}

\newtheorem{prop}{Proposition:}

\newtheorem{lem}{Lemma:}

\newtheorem{cor}{Corollary:}

\begin{document}
\large
\title{
{\bf 
Loewy structure of $G_1T$-Verma modules of singular highest weights
\footnotetext{\textit{2010 Mathematics Subject Classification.} 20G05.}
}
\thanks{supported in part by JSPS Grants in Aid for Scientific Research 
}
\author{
A\textsc{be} Noriyuki
\\
Hokkaido University
\\
Creative Research Institution (CRIS)
\\
abenori@math.sci.hokudai.ac.jp
\\
\and
K\textsc{aneda} Masaharu
\\
Osaka City University
\\
Department of Mathematics
\\
kaneda@sci.osaka-cu.ac.jp
}
}
\maketitle

\begin{abstract}
Let $G$ be a reductive algebraic group over an algebraically closed field of positive characteristic,
$G_1$ the Frobenius kernel of $G$, and $T$ a maximal torus of $G$.
We show that the parabolically induced
$G_1T$-Verma modules of singular highest weights are all rigid, determine their Loewy length, and describe their Loewy structure
using the periodic Kazhdan-Lusztig $Q$-polynomials.
We assume that the characteristic of the field is large enough that, in particular, Lusztig's conjecture for the irreducible $G_1T$-characters holds.

\end{abstract}

Let
$G$ be a reductive algebraic group over an algebraically closed field $\Bbbk$ of positive characteristic
$p$.
The Frobenius kernel $G_1$
of $G$ is an analogue of the Lie algebra of $G$ in characteristic $0$.
To keep track of weights, we 
consider
representations of $G_1T$
with $T$ a maximal torus of $G$.
In this paper we study $G_1T$-Verma modules, standard objects of the theory.


Many years ago Henning Andersen and the second author of the present paper showed that the $G_1T$-Verma modules of $p$-regular highest weights are all rigid of Loewy length
$1$
plus the dimension of the flag variety of $G$, and described their Loewy structure using  the periodic Kazhdan-Lusztig $Q$-polynomials
\cite{AK89}.
For that we assumed the validity of Lusztig's conjecture on the irreducible characters for $G_1T$-modules, or rather Vogan's equivalent version on the 
semisimplicity of certain $G_1T$-modules, modeling after Irving's method
\cite{I}, \cite{I88}.
Lusztig's conjecture is now a theorem for large $p$ as established by Andersen, Jantzen and Soergel
\cite{AJS}.
Pushing their graded representation theory, 
with a machinery of Beilinson, Ginzburg and Soergel
\cite{BGS}, 
we showed in \cite{AbK}
that the parabolic induction is graded on $p$-regular blocks, and 
determined the Loewy structure of parabolically induced
$G_1T$-Verma modules of
$p$-regular highest weights.
In this paper we use Riche's Koszulity of the $G_1$-block
algebras 
\cite{Ri}
to
uncover the structure of the 
parabolically induced
$G_1T$-Verma modules of $p$-singular highest weights, to complete the entire picture.


To describe our results precisely,
let us introduce some notations.
For simplicity we will assume throughout the paper that $G$ is simply connected and simple. 
Fix a Borel subgroup 
$B$
of $G$ containing $T$,
and 
choose a positive system $R^+$
of $R$ such that the roots of $B$ are
$-R^+$. Let $R^\rs$ denote the set of simple roots of $R^+$.
Let $\Lambda$ denote the weight lattice of $T$
equipped with a
partial order
such that
$\lambda\geq\mu$
iff
$\lambda-\mu\in\sum_{\alpha\in R^+}\bbN\alpha$.
Put
$\rho=\frac{1}{2}\sum_{\alpha\in R^+}\alpha$.
Let 
$W$
denote the Weyl group
of $G$ relative to $T$,
and let
$W_a=W\ltimes\bbZ R$ be the affine Weyl group 
with elements of $\bbZ R$ in $W_a$ acting on $\Lambda$ by translations.
We let $W_a$ 
act on
$\Lambda$ 
also via
$x\bullet\lambda=px\frac{1}{p}(\lambda+\rho)-\rho$,
$x\in W_a$, $\lambda\in \Lambda$.
Let
$R^\vee=\{\alpha^\vee\mid\alpha\in R\}$ denote the set of coroots of $R$,
and put $H_{\alpha,n}=\{
v\in\Lambda\otimes_\bbZ\bbR\mid\langle
v+\rho,\alpha^\vee\rangle=pn\}$,
$\alpha\in R$ and $n\in\bbZ$.
We call a connected component of
$(\Lambda\otimes_\bbZ\bbR)\setminus\cup_{\alpha\in R,n\in\bbZ}H_{\alpha,n}$
an alcove.
We say $\lambda\in\Lambda$ is $p$-regular iff it belongs to an alcove, otherwise 
$\lambda$ is $p$-singular.
Let also
$\Lambda^+=\{\lambda\in\Lambda\mid
\langle\lambda,\alpha^\vee\rangle\geq0\ \forall\alpha\in R^+\}$ the set of dominant weights.
We let $A^+=\{v\mid
\langle v+\rho,\alpha^\vee\rangle\in\mathopen]0,p\mathclose[\ \forall\alpha\in R^+\}$ denote the bottom dominant alcove.
For a closed subgroup $H$ of $G$ we let $H_1$ denote its Frobenius kernel.
Let $\hat\nabla=\ind_{B_1T}^{G_1T}$ denote the induction functor \cite[I.3]{J}
from the category of 
$B_1T$-modules to the category of 
$G_1T$-modules.
The $G_1T$-simple modules are parametrized by
their highest weights in $\Lambda$.
We let
$\hat L(\nu)$ denote the simple
$G_1T$-module of highest weight $\nu\in\Lambda$.
If $M$ is a 
finite dimensional 
$G_1T$-module,
we will write
$[M:\hat L(\nu)]$
for the composition factor multiplicity of
$\hat L(\nu)$ in $M$.

A Loewy filtration of
a 
finite dimensional 
$G_1T$-module
$M$ is a filtration of $M$ of minimal length
such that  each of its subquotients is semisimple.
The length of a Loewy filtration is uniform, called the Loewy length of $M$, 
denoted $\ell\ell(M)$.
Among the Loewy filtrations, the socle series of $M$ is defined by
 $0<\soc^1M<\soc^2M<\dots<\soc^{\ell\ell(M)}M=M$ with $\soc^1M=\soc M$, called the socle of $M$
 which is
the sum of simple submodules of $M$, and 
$\soc^iM/\soc^{i-1}M=\soc(M/\soc^{i-1}M)$ for
$i>1$.
Also the radical series of $M$ is defined by
 $0=\rad^{\ell\ell(M)}M<\dots<\rad^{2}M<\rad^1M<M$ with $\rad^1M=\rad M$, called the radical of $M$
 which is
the intersection of the maximal submodules of $M$,
and 
$\rad^iM=\rad(\rad^{i-1}M)$ for
$i>1$.
If $0<M^1<\dots<M^{\ell\ell(M)}=M$ is any Loewy filtration of $M$, 
$\rad^{\ell\ell(M)-i}M\leq M^i\leq\soc^iM$ for each $i$.
We say $M$ is rigid iff
the socle and the radical series of $M$ coincide.
We put $\soc_iM=\soc^iM/\soc^{i-1}M$
and
$\rad_jM=\rad^jM/\rad^{j+1}M$.

In this paper we show

\begin{thm}
Assume $p\gg0$.
Let $\nu\in\Lambda$ and let $N(\nu)$ denote the number of hyperplanes $H_{\alpha,n}$ on which $\nu$ lies.
The $G_1T$-Verma module $\hat\nabla(\nu)$ of highest weight $\nu$ is rigid of Loewy length
$1+\dim G/B-N(\nu)$.
If $x\in W_a$ is such that $\nu$ belongs to the upper closure of $x\bullet A^+$
and if
$\nu_0
=x^{-1}\bullet\nu$,
the Loewy structure of 
$\hat\nabla(\nu)$ is given by
\begin{multline*}
\sum_{i\in\bbN}
q^{\frac{\rd(y,x)-i}{2}}
[\soc_{i+1}\hat\nabla(
\nu):\hat L(y\bullet\nu_0)]
\\
=
\begin{cases}
Q^{y\bullet A^+,x\bullet A^+}
&\text{if $y\in W_a$ with $y\bullet\nu_0$
belonging to the upper closure of $y\bullet A^+$},
\\
0
&\text{else},
\end{cases}
\end{multline*}
where $\rd(y,x)$ is the distance from the alcove
$y\bullet A^+$ to the alcove
$x\bullet A^+$
\cite[1.4]{L80}
and
$Q^{y\bullet A^+,x\bullet A^+}$ is a polynomial from
\cite[1.8]{L80}.

\end{thm}

For this theorem to hold, we assume $p\gg0$ so that Lusztig's conjecture for the irreducible
characters of $G_1T$-modules and also the conditions
\cite[(10.1.1) and (10.2.1)]{Ri} from \cite{BMR}
hold.
While Fiebig \cite{F12} gives an explicit lower, as crude as it may be, bound of $p$
for Lusztig's conjecture to hold,
a recent work of Williamson
\cite{W} reveals that $p$ has, in general,  to be much bigger than $h$ the Coxeter number of $G$, which was the original bound 
for the conjecture to hold.
Compared to the restriction required for Lusztig's conjecture to hold, the other conditions in \cite{Ri}
are innocent.

We actually obtain, more generally, analogous results for parabolically induced module
$\hat\nabla_P(\hat L^P(\nu))=\ind_{P_1T}^{G_1T}(\hat L^P(\nu))$
with
$\hat L^P(\nu)$
denoting a simple $P_1T$-module
of highest weight $\nu$ for a parabolic subgroup $P$ of $G$.

For a category
$\cC$ we will denote the set of morphisms from object $X$ to $Y$ in $\cC$ by $\cC(X,Y)$ .

We are grateful to Wolfgang Soergel for suggesting us the problem and referring to \cite{Ri}.
Thanks are also due to Simon Riche for helpful comments to clarify the presentation of the paper.
The first author of the paper acknowledges the hospitality of Institute of Mathematics at Jussieu
and also Mittag-Leffler Institute, during the visit of which the work has been done.

\setcounter{equation}{0}
\begin{center}
$1^\circ$
{\bf 
Koszulity of the $G_1$-block algebras}

\end{center}

Throughout the paper we will
assume $p>h$ the Coxeter number of $G$ unless otherwise specified.
In particular, $p\Lambda\cap \bbZ R=p\bbZ R$.
All modules we consider are finite dimensional over $\Bbbk$.
Our basic strategy is to transport the known structure of a $G_1T$-block $\cC(\lambda)$ of
$p$-regular $\lambda\in\Lambda$ to
an arbitrary
block
$\cC(\mu)$
by the translation functor.
For $p\gg 0$, thanks to \cite{Ri},
the corresponding translation functor for the $G_1$-blocks is graded
and the $G_1$-block algebras 
are all Koszul.

\noindent
(1.1)
For $\nu\in\Lambda$
let $\hat L(\nu)$ denote the simple $G_1T$-module of highest weight
$\nu$, and 
$\hat P(\nu)$ 
the $G_1T$-projective cover of
$\hat L(\nu)$.
Let $\Omega$ be a $p$-regular orbit of
$W_a$ in $\Lambda$
and let
$\cC(\Omega)$ denote the corresponding
$G_1T$-block.
Thus $\cC(\Omega)=\cC(\nu)$, $\nu\in\Omega$, consists of $G_1T$-modules whose composition factors all have highest weights in
$
\Omega$.
Let 
$\Omega'$ be a system of representatives of 
$
\Omega$ under the translations by
$p\bbZ R$,
and let
$\hat P(\Omega)=\coprod_{\nu\in\Omega'}\hat P(\nu)$.
Then
$\coprod_{\gamma\in p\bbZ R}\cC(\Omega)(\hat P(\Omega)\otimes \gamma,
\hat P(\Omega))
$
forms a $p\bbZ R$-graded $\Bbbk$-algebra under the composition using the auto-functor $?\otimes\gamma$, $\gamma\in p\bbZ R$, on 
$\cC(\Omega)$.
If we let $\hat 
\bbE(\Omega)$ denote its opposite algebra,
$\coprod_{\gamma\in p\bbZ R}\cC(\Omega)(\hat P(\Omega)\otimes \gamma,
?)$ gives an equivalence of categories from 
$\cC(\Omega)$ to the category
of finite dimensional
$p\bbZ R$-graded $\hat
\bbE(\Omega)$-modules.
Moreover, $\hat
\bbE(\Omega)$ admits a $\bbZ$-grading compatibe with its $p\bbZ R$-gradation
\cite[18.17.1]{AJS}.
For $p$ large enough that Lusztig's conjecture holds, \cite[18.17]{AJS} has proved that 
$\hat
\bbE(\Omega)$
is Koszul with respect to its $\bbZ$-gradation.
Let us state Lusztig's conjecture in an equivalent form as follows:
$\forall x, y\in W_a$,
\[
\tag{L}
[\hat\nabla(x\bullet0):\hat L(y\bullet0)]
=
Q^{y\bullet A^+,x\bullet A^+}(1),
\]
where $Q^{y\bullet A^+,x\bullet A^+}$ is a polynomial from
\cite[1.8]{L80}.

Assuming (L), let $\tilde\cC(\Omega)$ denote the category of finite dimensional
$(p\bbZ R\times\bbZ)$-graded
$\hat
\bbE(\Omega)$-modules.
For each $\nu\in \Omega'$ let
$\tilde L(\nu)$ be the lift of $\hat L(\nu)$ in
$\tilde\cC(\Omega)$
as a direct summand of the degree 0 part of $\hat\bbE(\Omega)$.
If we denote 
the degree shift of objects in
$\tilde\cC(\Omega)$
by $[\gamma]$, $\gamma\in p\bbZ R$,
and by $\langle i\rangle$,
$i\in\bbZ$, any simple of $\tilde\cC(\Omega)$ may be written
$\tilde L(\nu)[\gamma]\langle i\rangle$,
$\nu\in\Omega'$, $\gamma\in p\bbZ R$, and $i\in\bbZ$,
in a unique way
up to isomorphism. 
As $\tilde L(\nu)[\gamma]$ is a lift of $\hat L(\nu+\gamma)=\hat L(\nu)\otimes\gamma$,
we will also write
$\tilde L(\nu+\gamma)$
for
$\tilde L(\nu)[\gamma]$.
For each 
$\nu\in\Omega$
the $G_1T$-Verma module
$\hat\nabla(\nu)$ of highest weight $\nu$ admits a lift
$\tilde\nabla(\nu)$ in $\tilde\cC(\Omega)$ such that its socle is
$\tilde L(\nu)$.
Likewise each projective
$\hat P(\nu)$ admits a lift
$\tilde P(\nu)$ which is the projective cover of
$\tilde L(\nu)$.

\noindent
(1.2)
Let
$\Lambda_p=\{
\nu\in\Lambda^+\mid
\langle\nu,\alpha^\vee\rangle<p\ \forall\alpha\in R^\rs\}$.
For 
$\nu\in\Lambda$ we write
$\nu=\nu^0+p\nu^1$
with $\nu^0\in\Lambda_p$ and
$\nu^1\in\Lambda$.
We let $L(\nu^0)$ denote the simple $G$-module of highest weight $\nu^0$,
which remains simple regarded as a $G_1$-module. 
All the simple $G_1$-modules are obtained thus.
One has
$\hat L(\nu)=L(\nu^0)\otimes p\nu^1$, 
$\hat P(\nu)=\hat P(\nu^0)\otimes p\nu^1$,
and $\hat P(\nu^0)$ provides the $G_1$-projective cover of
$L(\nu^0)$, which we will denote by
$P(\nu^0)$.

Let now
$\frg$ denote the Lie algebra of $G$,
$\bU\frg$ the universal enveloping algebra 
of $\frg$, 
and 
$(\bU\frg)_0$ 
the central reduction of $\bU\frg$
with respect to the Frobenius central character $0$.
As 
$(\bU\frg)_0$
coincides with the 
the algebra of distributions of $G_1$,
the representation theory of $G_1$ is equivalent to that of
$(\bU\frg)_0$.
For each $\nu\in\Lambda$ let
$(\bU\frg)_0^{\hat\nu}$ 
be the central reduction of
$(\bU\frg)_0$ with respect to the Harish-Chandra generalized character
$\hat\nu$.
This is the $G_1$-block component 
of $(\bU\frg)_0$
associated to
$\nu$.
Let $\cB(\nu)$ denote the category of finite dimensional
$(\bU\frg)_0^{\hat\nu}$-modules.
The algebra 
$(\bU\frg)_0^{\hat\nu}$ is equipped with a $\bbZ$-grading
\cite[6.3 and 10.2 line 16, p.\ 126]{Ri}.
We let 
$\cB^\gr(\nu)$
denote the category of finite dimensional graded
$(\bU\frg)_0^{\hat\nu}$-modules.
Let $\Lambda(\nu)=\{(w\bullet\nu)^0\mid w\in W\}$. 
Each $P(\eta)$, $\eta\in\Lambda(\nu)$,
admits a lift
$P^\gr(\eta)$ in
$\cB^\gr(\nu)$.
Let $P^\nu=\coprod_{\eta\in\Lambda(\nu)}P^\gr(\eta)$
and set
$\bbE(\nu)=
\cB(\nu)(P^\nu,P^\nu)^\op
$.
As $P^\nu$ is a projective
generator of $\cB(\nu)$
and as
$\bbE(\nu)= \coprod_{i\in\bbZ}
\cB^\gr(\nu)(P^\nu\langle i\rangle,P^\nu)
$
is equipped with a $\bbZ$-gradation, $\langle i\rangle$ denoting the degree shift,
$\cB(\nu)(P^\nu,\,?)$ induces 
an equivalence
from $\cB^\gr(\nu)$ to the category 
of finite dimensional $\bbZ$-graded
$\bbE(\nu)$-modules,
which we will denote by
$\tilde\cB(\nu)$.
For $p\gg0$, thanks to \cite[10.3]{Ri},
all $\bbE(\nu)$ are Koszul
by a careful choice of graded lift 
$P^\gr(\eta)$, $\eta\in\Lambda(\nu)$.

To be precise,
let 
$I\subseteq R^\rs$ and let $P$ denote the corresponding standard parabolic subgroup of
$G$ with the Weyl group
$W_I=\langle s_\alpha\mid\alpha\in I\rangle$,
where
$s_\alpha$ is the reflection associated to
$\alpha$.
Let
$\mu\in \Lambda$
lying in
the closure $\overline{A^+}$
of the alcove $A^+$
such that  $\rC_{W_a\bullet}(y\bullet\mu):=\{x\in W_a\mid xy\bullet\mu=y\bullet\mu\}=W_I$
for some
$y\in W_a$.
Let also $\lambda\in A^+$
such that
$\langle
y\bullet\lambda,\alpha^\vee\rangle=0$
$\forall\alpha\in I$. 
If $p\gg0$ so that 
the condition (L) holds,
one can take
each $P^\gr((w\bullet\lambda)^0)$ to satisfy
a certain condition \cite[8.1(\ddag)]{Ri}.
With this choice
\cite[Th. 9.5.1]{Ri} shows that
the graded algebra
$\bbE(\lambda)$ is Koszul.
For $\mu$
assume in addition to (L)
two more conditions, which go as follows:
the first one
\cite[10.1.1]{Ri}
coming from
\cite[Lem. 1.10.9(ii)]{BMR} reads, with $\cD^\lambda_{G/P}$
denoting the sheaf of PD-differential operators on $G/P$
twisted by the invertible sheaf
$\cL_{G/P}(\lambda)$,
\[
\tag{R1}
R^i\Gamma(G/P, \cD^\lambda_{G/P})=0
\quad\forall i>0.
\]
With $(\bU\frg)^\lambda$
denoting the central reduction of
$\bU\frg$
by 
the Harish-Chandra central 
character
$\lambda$, the second condition
\cite[10.2.1]{Ri}
coming also from
\cite[Lem. 1.10.9]{BMR}
reads that
\[
\tag{R2}
\text{
the natural morphism
$(\bU\frg)^\lambda\to\Gamma(G/P,\cD^\lambda_{G/P})$ be surjective.
}
\]
If $p\gg0$ so that (L), (R1) and (R2) all hold,
one can take each
$P^\gr(\eta)$, $\eta\in\Lambda(\mu)$, to satisfy
\cite[Th. 10.2.4]{Ri},
which makes $\bbE(\mu)$ also Koszul
\cite[Th. 10.3.1]{Ri}.
For any $\nu\in\Lambda$ there is 
by \cite[Lem. 1.5.2]{BMR}
some $\xi\in\Lambda$ such that
$\nu+p\xi\in W_a\bullet\mu$ 
with $\mu$ as above.
Thus under the conditions (L), (R1) and (R2)
we may assume 
that
all $G_1$-block algebras $\bbE(\nu)$ are Koszul.
For each $\eta\in\Lambda(\nu)$
we denote by
$\tilde L(\eta)$
the graded lift in $\tilde\cB(\nu)$
of $G_1$-simple $L(\eta)$
as a direct summand of $\bbE(\nu)_0$.
Let also
$\tilde P(\eta)=\coprod_{i\in\bbZ}\cB^\gr(\nu)(P^\nu\langle i\rangle, P^\gr(\eta))$ 
be a graded lift in
$\tilde\cB(\nu)$
of 
$P(\eta)$
to form the projective cover of $\tilde L(\eta)$.

\noindent
(1.3)
Assume from now on throughout the rest of the paper that
$p\gg0$ so that all the conditions
(L), (R1) and (R2) from (1.1) and (1.2) hold, unless otherwise specified.
Fix also $\lambda$ and $\mu$ as in (1.2).

For our purposes, as tensoring with $p\eta$, $\eta\in\Lambda$, is an 
equivalence from the $G_1T$-block
$\cC(\Gamma)$
of a $W_a$-orbit $\Gamma$
to the $G_1T$-block
$\cC(\Gamma+p\eta)$,
we have 
only to determine the structure of parabolically induced
$G_1T$-Verma modules
of highest weight
$x\bullet\mu$
with $\mu$ as above and $x\in
W_a$ such that $\langle x\rho,\alpha^\vee\rangle\in\mathopen]0,p\mathclose[$ $\forall\alpha\in R^\rs$.

If $\Omega=W_a\bullet\lambda$,
as $p>h$ by the standing hypothesis, $\bbE(\lambda)$ coincides 
by the linkage principle
with $\hat
\bbE(\Omega)$
from (1.1) as $\Bbbk$-algebras.
As two $\bbZ$-gradations on
the algebra must agree 
by their Koszulity
\cite[F.2]{AJS}, there is no ambiguity about the functor
from
$\tilde\cC(\Omega)$ to $\tilde\cB(\lambda)$ forgetting the
$p\bbZ R$-gradation,
which is compatible with the forgetful functor
from
the category of $G_1T$-modules to that of
$G_1$-modules.
Thus one has a commutative diagram of
forgetful functors
\[
\xymatrix{
\cC(\Omega)
\ar[d]
&&
\tilde\cC(\Omega)
\ar[d]
\ar[ll]
\\
\cB(\lambda)
&&
\tilde\cB(\lambda).
\ar[ll]
}
\]

\noindent
(1.4)
For each
$\nu\in\Lambda$
let
$\barpr_\nu$ denote the projection from the category 
of finite dimensional
$G_1$-modules to its
$\nu$-block
$\cB(\nu)$. 
For $\nu,\eta\in\overline{A^+}$
recall from \cite{BMR08}
the 
translation functor
$T_\nu^\eta=\barpr_\eta(L((\eta-\nu)^+)\otimes\ ?):
\cB(\nu)\to\cB(\eta)$ 
with
$(\eta-\nu)^+\in
W(\eta-\nu)\cap\Lambda^+$.

By \cite[Prop. 5.4.3 and Th. 6.3.4]{Ri} the adjoint 
translation functors 
$T_\lambda^\mu$ and $T_\mu^\lambda$ are graded
to form
a pair of functors
$
\cB^\gr(\lambda)\rightleftarrows
\cB^\gr
(\mu)$
such that
graded $T_\lambda^\mu$ is right
adjoint to graded $T_\mu^\lambda$. In turn, they induce a pair of graded 
functors,
which we will denote by 
$\tilde T_\lambda^\mu$ and $\tilde
T_\mu^\lambda$:
\begin{align*}
\tilde T_\mu^\lambda
&=
\coprod_{i\in\bbN}
\cB^\gr(\lambda)(P^\lambda\langle i\rangle, T_\mu^\lambda\,?)\circ
(P^\mu
\otimes_{\bbE(\mu)}?)
:
\tilde\cB(\mu)\to
\tilde\cB(\lambda),
\\
\tilde T^\mu_\lambda
&=
\coprod_{i\in\bbN}
\cB^\gr(\mu)(P^\mu\langle i\rangle, T^\mu_\lambda\,?)\circ
(P^\lambda
\otimes_{\bbE(\lambda)}?)
:\tilde\cB(\lambda)
\to
\tilde\cB(\mu)^\gr.
\end{align*}
Thus $\tilde T_\lambda^\mu$ is right adjoint to $\tilde T^\lambda_\mu$.

Let $N(\nu)$ denote the number of hyperplanes $H_{\alpha,n}$ on which $\nu\in\Lambda$ lies,
and put, in particular, $N
=N(\lambda)=\dim G/B$, $N_P
=N(\mu)=\dim G/P$.
A crucial fact to our results is Riche's 
\cite[10.2.8]{Ri} that asserts
for each
$w\in W$ with
$(w\bullet\mu)^0\in\Lambda(\mu)$, i.e., $(w\bullet\mu)^0$ belonging to the upper closure of an alcove containing 
some
$(w'\bullet\lambda)^0$, $w'\in W$,
$T_\mu^\lambda P^\gr((w\bullet\mu)^0)= P^\gr((w'\bullet\lambda)^0)\langle N-N_P\rangle$,
and hence
\begin{equation}
\tilde T_\mu^\lambda\tilde P((w\bullet\mu)^0)=\tilde P((w'\bullet\lambda)^0)\langle N-N_P\rangle.
\end{equation}

\setcounter{equation}{0}
\noindent
(1.5)
For each $\nu\in\Lambda$ 
let
$\widehat\pr_\nu$ denote the projection from the category 
of finite dimensional
$G_1T$-modules to the
block
$\cC(
\nu)$. 
For $\nu,\eta\in\overline{A^+}$
one has as in (1.4) the 
translation functor
$\hat T_\nu^\eta=\widehat\pr_\eta(L((\eta-\nu)^+)\otimes\ ?):
\cC(\nu)\to\cC(\eta)$ 
\cite[II.9.22]{J}.
Under the assumption $p>h$,
the functors
$\hat T_\lambda^\mu$ and $T_\lambda^\mu$ 
commute with the forgetful functors
as in \cite[Lem. 4.3.2]{Ri}:
\[
\xymatrix{
\cC(\lambda)
\ar[rr]^-{\hat T_\lambda^\mu}
\ar[d]
\ar
@{}
[drr]
|{\circlearrowright}
&&
\cC(\mu)
\ar[d]
\\
\cB(\lambda)
\ar[rr]_-{ T_\lambda^\mu}
&&
\cB(\mu).
}
\]

\setcounter{equation}{0}
\noindent
(1.6)
Under the forgetful functors,
$\hat\nabla=\ind_{B_1T}^{G_1T}$ yields an induction functor
$\bar\nabla=\ind_{B_1}^{G_1}$
from the category of 
$B_1$-modules to the category of 
$G_1$-modules.
If $M$ is a 
$G_1T$-module, it is semisimple iff it is semisimple as $G_1$-module
\cite[I.6.15]{J}.
Thus, in order to show that 
$\hat\nabla(x\bullet\mu)$, $x\in W_a$, is rigid,
we have only to show that $\bar\nabla(x\bullet\mu)$
is rigid.

For a facet
$F$ in $\Lambda\otimes_\bbZ\bbR$ with respect to
$W_a$
let
$\hat F$ denote 
its upper closure.
As
$\hat\nabla(x\bullet\mu)=\hat T_\lambda^\mu\hat\nabla(x\bullet\lambda)$,
$
\tilde T_\lambda^\mu\tilde\nabla(x\bullet\lambda)\in\tilde\cB(\mu)$ is
a graded lift of
$\bar\nabla(x\bullet\mu)$,
which we will denote by
$\tilde\nabla(x\bullet\mu)\langle i+N_P-N\rangle$
if
$x\bullet\mu\in\widehat{x'\bullet A^+}$,
$x'\in W_a$,
and if
$[\soc_{i+1}\hat\nabla(x\bullet\lambda):\hat L(x'\bullet\lambda)]=
[\tilde\nabla(x\bullet\lambda):\tilde L(x'\bullet\lambda)\langle i\rangle]\ne0$.
As $\bar\nabla(x\bullet\mu)$
has a simple socle and a simple head, so does its lift, and hence
the lift is rigid
by \cite[Prop. 2.4.1]{BGS}.
There now
follows
the rigidity of 
$\hat\nabla(x\bullet\mu)$.

\begin{prop}
All $G_1T$-Verma 
modules $\hat\nabla(
\nu)$, $\nu\in\Lambda$, are rigid.

\end{prop}

\setcounter{equation}{0}
\noindent
(1.7)
Let $w\in W$ and put ${^w\!B}=wBw^{-1}$,
$\hat\nabla_w=\ind_{{^w\!B}_1T}^{G_1T}$.
If $M$ is a 
$G_1T$-module,
let ${^w\!M}$ denote the $G_1T$-module $M$ with the $G_1T$-action twisted by $w$, i.e., we let
$g\in G_1T$ 
act on $m\in M$ by $w^{-1}gw$. 
For each $\nu\in\Lambda$
one has an isomorphism
$^w\hat\nabla(\nu)\simeq\hat\nabla_w(w\nu)$
\cite[II.9.3]{J}.
Thus

 \begin{cor}
All
$\hat\nabla_w(\nu)$, $w\in W$, $\nu\in\Lambda$, are rigid.
 
 \end{cor}

\setcounter{equation}{0}
\noindent
(1.8)
Let $J\subseteq R^\rs$, $Q$ the standard parabolic subgroup of $G$ associated to $J$
with the Weyl group denoted
$W_J$,
and let
$\hat\nabla_J=\ind_{Q_1T}^{G_1T}$ denote the induction functor from the category of
$Q_1T$-modules to the category of 
$G_1T$-modules.
Let $\nu\in\Lambda$
and let $\hat L^J(\nu)$ denote the simple $Q_1T$-module of highest weight $\nu$.
Choose a $p$-regular
$\eta\in\Lambda$ such that $\nu$ belongs to the upper closure of the $W_{J,a}$-alcove containing
$\eta$.
Under the Lusztig conjecture
(L)
we have shown in
\cite[3.9]{AbK} that
$\hat\nabla_J(\hat L^J(\eta))$ is graded, and in \cite[2.3]{AbK}
that
$\hat T_\eta^\nu(\hat\nabla_J(\hat L^J(\eta)))
\simeq
\hat\nabla_J(\hat L^J(\nu))$.
As 
$\hat\nabla_J(\hat L^J(\nu))$
has a simple head and socle
\cite[1.4]{AbK}, it follows again from
\cite[Prop. 2.4.1]{BGS} that

\begin{prop}
All parabolically induced $G_1T$-Verma modules
$\hat\nabla_J(\hat L^J(\nu))$,
$\nu\in\Lambda$,
are rigid.

\end{prop}

\setcounter{equation}{0}
\begin{center}
$2^\circ$
{\bf The Loewy structure}

\end{center}

Keep the notations from \S1.

\setcounter{equation}{0}
\noindent
(2.1)
For each $\nu\in \Lambda$ 
we will denote
$\hat L(\nu)$ by $\bar L(\nu)$ when 
regarded as a $G_1$-module.
Thus
$\bar L(\nu)=L(\nu^0)$.

\begin{lem}
Let $x\in W_a$.

(i)
One has
\[
\tilde T_\lambda^\mu\tilde L((x\bullet\lambda)^0)=
\begin{cases}
\tilde L((x\bullet\mu)^0)\langle N_P-N\rangle
&\text{if $x\bullet\mu\in
\widehat{x\bullet A^+}$, 
}
\\
0
&\text{else}.
\end{cases}
\]

(ii)
If $x\bullet\mu\in\widehat{x\bullet A^+}$, one has 
$
\forall i\in\bbN$,
\[
\hat T_\lambda^\mu\soc^i\hat\nabla(x\bullet\lambda)=
\soc^i\hat\nabla(x\bullet\mu).
\]

\end{lem}

\pf
(i)
We may by (1.5) assume that
$x\bullet\mu\in
\widehat{x\bullet A^+}$
 \cite[II.7.15, 9.22.4]{J},
 which occurs iff $(x\bullet\mu)^0$ lies in the upper closure of the alcove
$(x\bullet\lambda)^0$ belongs to.
Thus we are to show
in that case
that
$\tilde T_\lambda^\mu\tilde L((x\bullet\lambda)^0)=
\tilde L((x\bullet\mu)^0)\langle N_P-N\rangle
$.

As
$\tilde P((x\bullet\mu)^0)$
(resp. $\tilde P((x\bullet\lambda)^0)$) is a projective cover of
$\tilde L((x\bullet\mu)^0)$
(resp. $\tilde L((x\bullet\lambda)^0)$), we have for each
$n\in\bbZ$
\begin{align*}
\tilde\cB
(\mu)&
(\tilde P((x\bullet\mu)^0)\langle n\rangle, \tilde T_\lambda^\mu\tilde L((x\bullet\lambda)^0))
\simeq
\tilde\cB
(\lambda)(\tilde T^\lambda_\mu\tilde P((x\bullet\mu)^0)\langle n\rangle, \tilde L((x\bullet\lambda)^0))
\\
&\simeq
\tilde\cB
(\lambda)(\tilde P((x\bullet\lambda)^0)\langle n+N-N_P\rangle, \tilde L((x\bullet\lambda)^0))
\quad\text{by (1.4.1)},
\end{align*}
which is nonzero iff $n+N-N_P=0$,
and hence the assertion follows.

(ii)
Let
$\soc_{G_1}^i\bar\nabla((x\bullet\lambda)^0)$, $x\in W_a$, denote the $i$-th term of the $G_1$-socle series of
$\bar\nabla((x\bullet\lambda)^0)$,
which is just
$\soc^i\hat\nabla(x\bullet\lambda)$
regarded as $G_1$-module.
As the socle series and the gradation over $\bbE(\lambda)$
(resp. $\bbE(\mu)$) coincide on
$\tilde\nabla((x\bullet\lambda)^0)$
(resp. 
$\tilde\nabla((x\bullet\mu)^0)$)
by \cite[Prop. 2.4.1]{BGS}, we see from (i)  that
$T_\lambda^\mu\soc_{G_1}^i\bar\nabla((x\bullet\lambda)^0)=
\soc_{G_1}^i\bar\nabla((x\bullet\mu)^0)$, and hence
the assertion.

 \setcounter{equation}{0}
\noindent
(2.2)
$\forall x, y\in W_a$, let
$Q^{y\bullet
A^+, x\bullet A^+}(q)=\sum_jQ^{y,x}_jq^{\frac{j}{2}}\in\bbZ[q]$
be the periodic 
Kazhdan-Lusztig 
$Q$-polynomial from \cite{L80}.
Put
$Q^{y,x}=Q^{y\bullet
A^+, x\bullet A^+}(q)$
for simplicity.
Recall from
\cite{AK89},
\cite[18.19]{AJS}/\cite[5.1, 2]{AbK}
\[
\sum_{i\in\bbN}
q^{\frac{\rd(y,x)-i}{2}}
[\soc_{i+1}\hat\nabla(x\bullet\lambda):\hat L(y\bullet\lambda)]
=
\sum_{i\in\bbN}
q^{\frac{\rd(y,x)-i}{2}}
[
\tilde\nabla(x\bullet\lambda):\tilde L(y\bullet\lambda)\langle-i\rangle]
=
Q^{y,x},
\]
where
$\rd(y,x)=\rd(y\bullet A^+,x\bullet A^+)$
is the distance from the alcove
$y\bullet A^+$ to the alcove
$x\bullet A^+$
\cite{L80}.
Let
$W_a(\mu)=\{x\in W_a\mid x\bullet\mu\in\widehat{x\bullet A^+}\}$.
For each
$x\in W_a(\mu)
$,
(2.1.ii) shows that
\[
\sum_{i\in\bbN}
q^{\frac{\rd(y,x)-i}{2}}
[\soc_{i+1}\hat\nabla(x\bullet\mu):\hat L(y\bullet\mu)]
=
\begin{cases}
Q^{y,x}
&\text{if $y\in W_a(\mu)$},
\\
0
&\text{else}.
\end{cases}
\]

\setcounter{equation}{0}
\noindent
(2.3)
One can do the same with parabolically induced $G_1T$-Verma modules 
$\hat\nabla_J(\hat L^J(\nu))$,
$J\subseteq R^\rs$,
$\nu\in\Lambda$,
from (1.8),
via
\cite[2.3]{AbK}.
Let
$W_{J,a}=W_J\ltimes\bbZ R_J$
denote the affine Weyl group
for $P_J$.

\begin{thm}
Let
$\nu\in\Lambda$,
$x\in W_a$ such that $\nu\in\widehat{x\bullet A^+}$,
and put
$\nu_0
=x^{-1}\bullet\nu$.
Then
\begin{multline*}
\sum_{i\in\bbN}
q^{\frac{\rd(y,x)-i}{2}}
[\soc_{i+1}\hat\nabla_J(
\hat L^J(\nu)):\hat L(y\bullet\nu_0)]
\\
=
\begin{cases}
\sum_{z\in W_{J,a}}Q^{y\bullet A^+,z\bullet A^+}
(-1)^{\rd_J(z,x)}\hat
P^J_{z\bullet A^+,x\bullet A^+}
&\text{if
$y\in
W_a(\mu)$},
\\
0
&\text{else},
\end{cases}
\end{multline*}
where $\hat{P}^J_{z\bullet A^+,x\bullet A^+}$ is a $\hat P$-polynomial from
\cite{Kato}
for $W_{J,a}$ and $\rd_J(z,x)$ is the distance from
$z\bullet A^+$ to
$x\bullet A^+$ with respect to
$W_{J,a}$.

\end{thm}

\setcounter{equation}{0}
\noindent
(2.4)
Finally we determine the Loewy length
of all parabolically induced $G_1T$-Verma  modules.
We first need analogues of
\cite[Props. 3.2 and 3.3]{I}.

Let
$\hat\Delta(\nu)=\Dist(G_1)\otimes_{\Dist(B_1)}\nu=\hat\nabla(\nu)^\tau$
the $\Bbbk$-linear dual of
$\hat\nabla(\nu)$ twisted by the Chevalley anti-involution $\tau$ of $G$
\cite[II.2.12]{J}.
We say a $G_1T$-module $M$ admits a $\hat\nabla$-filtration
iff there is a filtration
$0=M^0<M^1<\dots<M^r=M$
of $G_1T$-modules with
each $M^i/M^{i-1}\simeq\hat\nabla(\nu_i)$  for some
$\nu_i\in\Lambda$,
in which case
one can arrange the filtration such that
$\nu_i\not<\nu_j$ if $i>j$
\cite[II.9.8]{J}.
Whenever $M$ admits a $\hat\nabla$-filtration, we will assume that such a rearrangement has been done.

Let
$W_\nu=\rC_{W_a\bullet}(\nu)$
and take an alcove $A$ in the closure of which $\nu$ lies.
Choose $\eta\in \Lambda$ in $A$.
Let
$\eta^+$
(resp. $\eta^-$)
denote the highest (resp. lowest)
weight in
$W_\nu\bullet\eta$.
Let us also denote by
$\hat T_\nu^\eta:\cC(W_a\bullet\nu)\to\cC(W_a\bullet\eta)$
and $\hat T^\nu_\eta:\cC(W_a\bullet\eta)\to\cC(W_a\bullet\nu)$
the 
associated translation functors.

\begin{lem}
Assume
$p\gg0$ so that
(L) holds.

(i)
$\hat\Delta(\eta^+
)\leq
\rad^{N(\nu)
}
\hat T^\eta_\nu\hat\Delta(\nu
)$.

(ii)
$\hat L(\eta^-)
\leq\soc_{N(\nu)+1}\hat\nabla(\eta^+)$.

(iii)
$\ell\ell(\hat T_\nu^\eta\hat L(\nu
))\geq
2N(\nu)+1$.

(iv)
$\forall M\in\cC(
\nu)$,
$\ell\ell(\hat T_\nu^\eta
M)\geq
2N(\nu)+\ell\ell(M)$.

\end{lem}

\pf
(i)
Recall from
\cite[18.13]
{AJS} that
the translation functors
$\hat T_\eta^\nu$ and $\hat T_\nu^\eta$
admit graded versions, 
denoted 
$T_!$ and $T^*$, resp.
If we let
$\tilde\Delta(\eta)$
denote the graded version of
$\hat\Delta(\eta)$,
$T^*T_!\tilde\Delta(\eta^-)$
admits by
\cite[18.15]{AJS}
a filtration
with the subquotients
$\tilde\Delta(w\bullet\eta^-)\langle
o(w\bullet\eta^-)\rangle$,
$w\in W_\nu$,
where
$o(w\bullet\eta^-)$ denotes the number of
hyperplanes
$H_{\alpha,n},\alpha\in R^+,n\in\bbZ$, on which $\nu$ lies and such that
$w\bullet\eta^-$ belongs to their positive sides
\cite[15.13]{AJS}.
Thus the graded version
of $\hat L(\eta^+)=\hd\hat\Delta(\eta^+)$
appears in
$T^*T_!\tilde\Delta(\eta^-)$
as
$
\tilde L(\eta^+)\langle N(\nu)\rangle$
while that
of
$\hat L(\eta^-)=\hd\hat\Delta(\eta^-)=
\hd\hat T_\nu^\eta\hat\Delta(\nu)$
appears as
$\tilde L(\eta^-)$.
Under the assumtion (L),
$\tilde\Delta(\eta^-)$
is 
graded over the Koszul algebra
$\hat\bbE(W_a\bullet\eta)$ from (1.1),
and
so therefore is 
$T^*T_!\tilde\Delta(\eta^-)$.
As $\hat T_\nu^\eta\hat\Delta(\nu)$ has a simple socle and a simple head, 
its Loewy series coincides with the grading filtration up to degree shift by
\cite{BGS}. 
It follows 
that
$\hat L(\eta^+)$ appears in $\rad_{N(\nu)}\hat T_\nu^\eta\hat\Delta(\nu)$,
and hence
$\hat\Delta(\eta^+)\leq\rad^{N(\nu)}\hat T_\nu^\eta\hat\Delta(\nu)$.

(ii)
Note first that the number of times
$\hat\nabla(\eta^+
)$ appears in the $\hat\nabla$-filtration of 
$\hat P(\eta^-
)$
is by the translation principle equal to
\begin{equation}
[\hat\nabla(\eta^+
)
:
\hat L(\eta^-
)]=1.
\end{equation}
Thus, if 
$r=\max\{i\in\bbN\mid\eta_i=
\eta^+
\}$
in the $\hat\nabla$-filtration 
$M^\bullet$
of
$\hat P(\eta^-)$
with the suquotients
$M^i/M^{i-1}\simeq
\hat\nabla(\eta_i)$,
one has
$\hat T_\nu^\eta\hat\nabla(
\nu
)\leq M^r$.
By (1) and by
\cite[3.5]{AK89}
there is unique $ j\in\bbN
$ such that
$[\soc_{j+1}
M^r
:\hat L(\eta^+
)]
=
[\soc_{j+1}\hat\nabla(\eta^+
):\hat L(\eta^-
)]=1$.
As
$[\hat T_\nu^\eta\hat\nabla(\nu
):\hat L(\eta^+
)]\ne0$,
we must have
$[\soc_{j+1}\hat T_\nu^\eta\hat\nabla(\nu
):\hat L(\eta^+
)]=1$.
Then 
taking the $\tau$-dual yields that
$[\rad_{j
}\hat T_\nu^\eta\hat\Delta(\nu
):\hat L(\eta^+
)]=1$,
and hence
$j=N(\nu)$
by (i).

(iii)
Consider a filtration of
$\hat T_\nu^\eta\hat\Delta(
\nu
)$
with the subquotients
$\hat\Delta(w\bullet\eta)$,
$w\in W_\nu$.
By the weight consideration
$\hat T_\nu^\eta\hat L(
\nu
)$ must contain all the composition factors of
$\hat T_\nu^\eta\hat\Delta(
\nu
)$
isomorphic to $\hat L(\eta^-
)$.

On the other hand,
$[\hat\Delta(w\bullet\eta^+):\hat L(\eta^-
)]=1$
$\forall w\in W_\nu$
as in (1).
Thus $\hat T_\nu^\eta\hat L(\nu
)$ contains a composition factor
$\hat L(\eta^-)$ corresponding to one in each of
$\hat \Delta(w\bullet\eta^+)$,
$w\in W_\nu$.
Consider 
the factor corresponding to the one in
$\hat\Delta(\eta^+)$.
Let $\theta\in
G_1T\Mod(\hat\Delta(\eta^+
),
\hat T_\nu^\eta \hat L(\nu
))$
be the restriction to
$\hat\Delta(\eta^+
)$ of the quotient
$\hat T_\nu^\eta \hat \Delta(\nu
)
\to
\hat T_\nu^\eta \hat L(\nu
)$.
Then
$\im\theta\leq\rad^{N(\nu)
}\hat T_\nu^\eta
\hat L(\nu
)$ by (i).
As the composition factor $\hat L(\eta^-
)$ comes from the one in
$\rad_{N(\nu)
}\hat\Delta(\eta^+
)$
by (ii),
it lies in
$\rad_{N(\nu)
}(\im\theta)$.
It follows that
$2N(\nu)+1\leq\ell\ell(\im\theta)+N(\nu)
\leq
\ell\ell(\hat T_\nu^\eta
\hat L(\nu
))$.

(iv)
Consider a nonsplit exact sequence
$0\to\hat L(y\bullet\nu)\to
K\to\hat L(x\bullet\nu)\to0$,
$x, y\in W_a$,
with
$x\bullet\nu>y\bullet\nu$.
There is 
an epi
$\hat\Delta(x\bullet\nu)\twoheadrightarrow K$.
As
$\hat T_\nu^\eta\hat\Delta(x\bullet\nu)$ has a simple head, so does
$\hat T_\nu^\eta
K$.
In particular,
$\hat T_\nu^\eta K$ is indecomposable, and so therefore is
$(\hat T_\nu^\eta K)^\tau\simeq\hat T_\nu^\eta(K^\tau)$.

We now argue by induction on
$\ell\ell(M)$.
We may assume $M$ has a simple head.
Let
$\hat L(x\bullet\nu)=\hd M$,
$x\in W_a$.
Take a quotient $M/M'$ with
$\rad M>M'>\rad^2M$
which fits in a short exact sequence
$0\to\hat L(y\bullet\nu)\to
M/M'\to\hat L(x\bullet\nu)\to0$
for some
$y\in W_a$.
As
$\hat T_\nu^\eta(M/M')$ is indecomposable,
the exact sequence
$0\to\hat T_\nu^\eta(\rad M)\to
\hat T_\nu^\eta M\to\hat T_\nu^\eta\hat L(x\bullet\nu)\to0$ cannot split.
Thus
$\ell\ell(\hat T_\nu^\eta M)\geq
\ell\ell(\hat T_\nu^\eta(\rad M))+1$, as desired.

\setcounter{equation}{0}
\noindent
(2.5)
Keep the notation of (2.4).
Let
$w_0$ denote the longest element of
$W$.
$\forall x\in W_a$,
recall from
\cite[3.4.2]{AK89}
that
\begin{equation}
\ell\ell\hat P(\nu
)\geq
2\ell\ell(\hat\nabla(w_0\bullet\nu))-1
\geq
2N-2N(\nu)+1,
\end{equation}
and from
\cite[2.3]{AK89}
that
$\ell\ell\hat\nabla(w_0\bullet\nu)\geq
N-N(\nu)+1$.
Thus
\begin{align*}
2N+1
&=
\ell\ell\hat P(\eta^-
)
\quad\text{by \cite[5.4]{AK89}}
\\
&=
\ell\ell(\hat T_\nu^\eta\hat P(\nu
))
\geq
\ell\ell(\hat P(\nu
))+2N(\nu)
\quad\text{by (2.4)}
\\
&\geq
2N-2N(\nu)+1+2N(\nu)
=
2N+1.
\end{align*}
It follows that
$\ell\ell\hat P(\nu
)=2N-2N(\nu)+1$,
and then
$\ell\ell\hat\nabla(w_0\bullet\nu)=N-N(\nu)+1$ by (1).
As
$\hat\nabla_w((w\bullet\nu)\langle w\rangle)
\simeq
{^w\hat\nabla}(\nu)\otimes
p(w\bullet0)$
$\forall w\in W$
by (1.5),
we have
\[
\ell\ell\hat\nabla_w((w\bullet\nu)\langle w\rangle)
=1+N-N(\nu)=
1+\dim G/B-N(\nu).
\] 
Let us also record

\begin{thm}
Assume
$p\gg0$ so that
(L) holds.
$\forall\nu\in\Lambda$,
$\ell\ell(\hat P(\nu))=2N-2N(\nu)+1$.

\end{thm}

\setcounter{equation}{0}
\noindent
(2.6)
Recall the notation of (1.8).
To find the Loewy length of
$\hat\nabla_J(\hat L^J(\nu))$,
we first recall some identities from
\cite{AbK}.
These hold without restrictions on $p$.
Let $w_J$ denote the longest element of
$W_J$
and put
$w^J=w_0w_J$.
Let
$\nu\in\Lambda$.
We will write
$\nu\langle w\rangle$,
$w\in W$, for
$\nu+(p-1)(w\bullet 0)$.

One can reformulate
\cite[1.4]{AbK} as an isomorphism
$\hd\hat\nabla_J(\hat L^J(\nu))
\simeq
L((w^J\bullet\nu)^0)\otimes_\Bbbk
p\{(w^J)^{-1}\bullet(w^J\bullet\nu)^1\}$.
Also, \cite[4.5]{AbK} carries over to arbitrary $\nu\in\Lambda$:
$\hd{^{w^J}\hat\nabla}_J(\hat L^J(\nu))
\simeq
\hat
L(w^J\bullet\nu)\otimes_\Bbbk
\{-p(w^J\bullet0)\}$.
We then have a commutative diagram from
\cite[4.6.1]{AbK}
\[
\xymatrix{
\hat\nabla_{w^J}((w^J\bullet
\nu)\langle
w^J\rangle)\otimes\{-p(w^J\bullet0)\}
\ar[rrr]^-{\phi_{w^J}\otimes\{-p(w^J\bullet0)\}
}
&&&
\hat\nabla(w^J\bullet
\nu)\otimes\{-p(w^J\bullet0)\}
\\
^{w^J}\hat\nabla_J(\hat L^J(\nu))
\ar@{->>}[rrr]
\ar@{^(->}!<0ex,2.5ex>;[u]
&&&
\hat
L(w^J\bullet
\nu)\otimes\{-p(w^J\bullet0)\}
\ar@{^(->}[u]
}
\]
and another from \cite[4.6.3]{AbK}
\[
\xymatrix{
\hat\nabla_{w_0}((w^J\bullet
\nu)\langle
w_0\rangle)\otimes\{-p(w^J\bullet0)\}
\ar[dd]_-{\phi'_{w^I}\otimes\{-p(w^J\bullet0)\}
}
\ar@{->>}[rd]
\\
&
^{w^J}\hat\nabla_J(\hat L^J(\nu))
\ar@{_(->}!<2ex,0ex>;[ld]
\\
\hat\nabla_{w^J}((w^J\bullet\nu)\langle
w^J\rangle)\otimes\{-p(w^J\bullet0)\}.
}
\]
If we write
$w^J=s_{i_1}\dots
s_{i_n}$ in a reduced expression with
$s_i$ denoting the reflection associated to the simple
root $\alpha_i$,
the homomorphism
$\phi_{w^J}:\hat\nabla_{w^J}((w^J\bullet
\nu)\langle
w^J\rangle)\to
\hat\nabla(w^J\bullet
\nu)$ is the composite
\begin{multline*}
\hat\nabla_{s_{i_1}\dotsm
s_{i_n}}((w^J\bullet
\nu)\langle
s_{i_1}\dotsm
s_{i_n}\rangle)\to
\hat\nabla_{s_{i_1}\dotsm
s_{i_{n-1}}}((w^J\bullet
\nu)\langle
s_{i_1}\dotsm
s_{i_{n-1}}\rangle)
\to\dotsb
\\\to
\hat\nabla_{s_{i_1}
s_{i_2}}((w^J\bullet
\nu)\langle
s_{i_1}
s_{i_2}\rangle)\to
\hat\nabla_{s_{i_1}}((w^J\bullet
\nu)\langle
s_{i_1}\rangle)
\to
\hat\nabla((w^J\bullet
\nu))
\end{multline*}
with each
$\hat\nabla_{s_{i_1}\dotsm
s_{i_r}}((w^J\bullet
\nu)\langle
s_{i_1}\dotsm
s_{i_r}\rangle)\to
\hat\nabla_{s_{i_1}\dotsm
s_{i_{r-1}}}((w^J\bullet
\nu)\langle
s_{i_1}\dotsm
s_{i_{r-1}}\rangle)
$
bijective iff 
$\langle
w^J\bullet\nu+\rho,s_{i_1}\dotsm
s_{i_{r-1}}\alpha_{i_r}^\vee\rangle\equiv0\mod p$
\cite[2.2]{AK89}.
Thus, if we put
$R^+(w)=\{\alpha\in R^+\mid w
\alpha<0\}$,
$w\in W$, and
$R^+_\nu=\{\alpha\in R^+\mid 
\langle\nu+\rho,\alpha^\vee\rangle\equiv0\mod p\}$,
then
$\phi_{w^J}\otimes\{-p(w^J\bullet0)\}$
annihilates
$\soc^{\ell(w^J)-|R^+(w^J)\cap
R^+_\nu|}\hat\nabla_{w^J}((w^J\bullet
\nu)\langle
w^J\rangle)\otimes\{-p(w^J\bullet0)\}
$,
and hence
\begin{equation}
\ell\ell
^{w^J}\hat\nabla_J(\hat L^J(\nu))
\geq
\ell(w^J)-|R^+(w^J)\cap
R^+_\nu|+1.
\end{equation}
Likewise,
$\phi'_{w^J}\otimes\{-p(w^J\bullet0)\}$
annihilates
$\soc^{\ell(w_J)-|(R^+\setminus
R^+(w^J))\cap
R^+_\nu|}\hat\nabla_{w_0}((w^J\bullet
\nu)\langle
w_0\rangle)\otimes\{-p(w^J\bullet0)\}
$.

We now 
assume (L) again.
As $\ell\ell\hat\nabla_{w_0}((w^J\bullet
\nu)\langle
w_0\rangle)
=
N-N(\nu)+1$
by (2.5),
\begin{equation}
\ell\ell
^{w^J}\hat\nabla_J(\hat L^J(\nu))
\leq
N-N(\nu)+1-
\{\ell(w_J)-|(R^+\setminus
R^+(w^J))\cap
R^+_\nu|\}.
\end{equation}
As $N(\nu)=|R^+(w^J)\cap
R^+_\nu|+|(R^+\setminus
R^+(w^J))\cap
R^+_\nu|$, it now follows from
(1) and (2) that
\begin{align*}
\ell\ell
\hat\nabla_J(\hat L^J(\nu))
&=\ell(w^J)-|R^+(w^J)\cap
R^+_\nu|+1
=
1+\dim G/Q
-|R^+(w^J)\cap
R^+_\nu|
\\
&
=|R^+(w^J)|-|R^+(w^J)\cap
R^+_\nu|+1
=
1+|R^+(w^J)\setminus
R^+_\nu|.
\end{align*}
Thus

\begin{thm}
Assume
$p\gg0$ so that
(L) holds.
All $\hat\nabla_J(\hat L^J(\nu))$,
$J\subseteq R^\rs$, $\nu\in\Lambda$, have Loewy length
$1+
|R^+(w^J)\setminus
R^+_\nu|
$.

\end{thm}


\begin{thebibliography}{AAAAA}

\bibitem[AbK]{AbK}  
Abe, N. and Kaneda, M., 
{\it
On the structure of parabolically induced $G_1T$-Verma modules}, JIM Jussieu
{\bf 14} Issue 01 (2015), 185-220


\bibitem[AJS]{AJS} Andersen, H.H., Jantzen, J.C. and Soergel, W., Representations of quantum groups at a $p$-th root of unity and of semisimple groups in characteristic $p$ : independence of $p$, Ast\' erisque {\bf 220}, 1994 (SMF)






\bibitem[AK]{AK89}  
Andersen, H.H. and Kaneda M., 
{\it
Loewy series of modules for the first Frobenius kernel in a reductive algebraic group},
Proc. LMS (3) {\bf 59} (1989), 
74--98

























\bibitem[BGS]{BGS} Beilinson, A., Ginzburg, D. and Soergel, W., {\it
Koszul Duality Patterns in Representation Theory},
J. AMS
{\bf 9} No.2 (1996),
473-527








\bibitem[BMR06]{BMR}  
Bezrukavnikov, R.,
Mirkovic, I.
and Rumynin, D.,
{\it Singular localization and intertwining functors  for reductive Lie algebras in prime characteristic}, 
Nagoya Math. J.
{\bf 184} (2006), 1--55


\bibitem[BMR08]{BMR08}  Bezrukavnikov, R., Mirkovic, I. and Rumynin, D., {\it Localization of modules for a semisimple Lie algebra in prime characteristic}, Ann. Math. {\bf 167} (2008), 945--991

































\bibitem[F]{F12} P. Fiebig, {\it An upper bound on the exceptional characteristics for Lusztig's character formula},  J. reine angew. Math.
{\bf 673} (2012), 1-31





















































\bibitem[I85]{I} Irving, R. S., {\it Projective modules in the category
$\cO_S$: Loewy series}, Trans. AMS
{\bf 291} No. 2 (1985), 733--754

\bibitem[I88]{I88} Irving, R. S., {\it The socle filtration of a Verma module}, Ann. Sci. ENS
{\bf 21} (1988), 47--65


\bibitem[J]{J}
Jantzen, J. C., Representations of Algebraic Groups, 2003
(American Math. Soc.)
















































\bibitem[Kat]{Kato} Kato S., {\it On the Kazhdan-Lusztig polynomials for affine Weyl groups}, Adv. Math. {\bf 55} (1995), 103-130



\bibitem[L80]{L80}
Lusztig, G.,
{\it Hecke algebras and Jantzen's generic decomposition patterns},
Adv. Math.
{\bf 37} (1980),
121-164







\bibitem[Ri]{Ri} Riche, S., {\it Koszul duality and modular representations of semisimple Lie algebras}, Duke Math. J. {\bf 154}
(2010), 31-134


















\bibitem[W]{W} Williamson, G., {\it Schubert calculus and torsion}, arXiv:1309.5055



























































































































\end{thebibliography}
\end{document}